\documentclass[a4paper]{article}

\usepackage[margin=1in]{geometry}

\usepackage[T1]{fontenc}
\newcommand{\changefont}[3]{
\fontfamily{#1} \fontseries{#2} \fontshape{#3} \selectfont}

\changefont{ptm}{m}{n}

\usepackage{setspace} 
 \usepackage{graphicx}    % needed for including graphics e.g. EPS, PS

\newcommand \be{\begin{equation}}
\newcommand \ee{\end{equation}}
\newcommand \ba{\begin{eqnarray}}
\newcommand \ea{\end{eqnarray}}

\def\bit{\begin{itemize}}
\def\eit{\end{itemize}}

\usepackage{amsfonts}
\usepackage{amssymb}
\usepackage{graphics}
\usepackage{mathrsfs}
\usepackage{color}
\newtheorem{remark}{Remark}[section]

\newtheorem{theorem}{Theorem}[section]

\newtheorem{lemma}{Lemma}[section]

\newtheorem{definition}{Definition}[section]

\long\def\symbolfootnote[#1]#2{\begingroup%
\def\thefootnote{\fnsymbol{footnote}}\footnote[#1]{#2}\endgroup} 

\begin{document}

%\begin{frontmatter}
%

\begin{center}
\Large \textbf{Unpredictable Points and Chaos}
\end{center}

\begin{center}
\normalsize \textbf{Marat Akhmet$^{a,}\symbolfootnote[1]{Corresponding Author Tel.: +90 312 210 5355,  Fax: +90 312 210 2972, E-mail: marat@metu.edu.tr}$ and Mehmet Onur Fen$^a$} \\
\vspace{0.2cm}
\textit{\textbf{\footnotesize $^a$Department of Mathematics, Middle East Technical University, 06800, Ankara, Turkey}} 
\vspace{0.1cm}
\vspace{0.1cm}
\end{center}

\vspace{0.3cm}

\begin{center}
\textbf{Abstract}
\end{center}

\noindent\ignorespaces

It is revealed that a special kind of Poisson stable point, which we call an unpredictable point, gives rise to the existence of chaos in the quasi-minimal set. The existing definitions of chaos are formulated in sets of motions. This is the first time that description of chaos is initiated from a single motion. The theoretical results are exemplified by means of the symbolic dynamics.

\vspace{0.2cm}
 
\noindent\ignorespaces \textbf{Keywords:} Unpredictable point; Chaos; Poisson stability; Symbolic dynamics; Quasi-minimal set

\vspace{0.6cm}

\section{Introduction}

The mathematical dynamics theory, which was founded by Poincar\'e \cite{poincare} and significantly developed by the french genius and Birkhoff \cite{birkhoff}, was a source as well as the basis for the later discoveries and thorough investigations of complex dynamics \cite{henon}-\cite{smale}. The homoclinic chaos was discussed by Poincar\'e \cite{poincare1}, and Lorenz \cite{lorenz} observed that a strange attractor contains a Poisson stable trajectory. Possibly, it was Hilmy \cite{Hilmy36,Nemytskii} who gave the first definition of a quasi-minimal set (as the closure of the hull of a Poisson stable motion). In \cite{Nemytskii} one can find a theorem by Hilmy, which states the existence of an uncountable set of Poisson stable trajectories in a quasi-minimal set. We modify the Poisson stable points to unpredictable points such the quasi-minimal set is chaotic.  

Let $(X, d)$ be a metric space and $\mathbb T$ refer to either the set of real numbers or the set of integers. A mapping $f: \mathbb T \times X \to X$ is a flow on $X$ \cite{sell} if:
\begin{itemize}
\item[(i)] $f(0,p)=p$ for all $p \in X;$
\item[(ii)] $f(t,p)$ is continuous in the pair of variables $t$ and $p;$
\item[(iii)] $f(t_1, f(t_2,p))=f(t_1+t_2,p)$ for all $t_1,$ $t_2 \in \mathbb T$ and $p \in X.$ 
\end{itemize}  
If a mapping $f: \mathbb T_+ \times X \to X,$ where $\mathbb T_+$ is either the set of non-negative real numbers or the set of non-negative integers, satisfies $(i),$ $(ii)$ and $(iii),$ then it is called a semi-flow on $X$ \cite{sell}. 
 
Suppose that $f$ is a flow on $X.$ A point $p \in X$ is stable $P^+$ (positively Poisson stable) if for any neighborhood $\mathcal{U}$ of $p$ and for any $H_1>0$ there exists $t\ge H_1$ such that $f(t,p) \in \mathcal{U}.$  Similarly, $p \in X$ is stable $P^-$ (negatively Poisson stable) if for any neighborhood $\mathcal{U}$ of $p$ and for any $H_2<0$ there exists $t\le H_2$ such that $f(t,p) \in \mathcal{U}.$ A point $p \in X$ is called stable $P$ (Poisson stable) if it is both stable $P^+$ and stable $P^-$ \cite{Nemytskii}.

For a fixed $p \in X,$ let us denote by $\Omega_p$ the closure of the trajectory $\mathcal{T}(p) = \left\{ f(t,p): t \in \mathbb T \right\},$ i.e., $\Omega_p = \overline{\mathcal{T}(p)}.$ The set $\Omega_p$ is a quasi-minimal set if the point $p$ is stable P and $\mathcal{T}(p)$ is contained in a compact subset of $X$ \cite{Nemytskii}. We will also denote $\Omega^+_p = \overline{\mathcal{T}^+(p)},$ where $\mathcal{T}^+(p) = \left\{ f(t,p): t \in \mathbb T_+ \right\}$ is the positive semi-trajectory through $p.$  
 
An essential result concerning quasi-minimal sets was given by Hilmy \cite{Hilmy36}. It was demonstrated that if the trajectory corresponding to a Poisson stable point $p$ is contained in a compact subset of $X$ and $\Omega_p$ is neither a rest point nor a cycle, then $\Omega_p$ contains an uncountable set of motions everywhere dense and Poisson stable. The following theorem can be proved by adapting the technique given in \cite{Hilmy36,Nemytskii}. 

\begin{theorem} \label{thm1a}   
Suppose that $p \in X$ is stable $P^+$ and $\mathcal{T}^+(p)$ is contained in a compact subset of $X.$ If \ $\Omega^+_p$ is neither a rest point nor a cycle, then it contains an uncountable set of motions everywhere dense and stable $P^+.$
\end{theorem}

\section{Unpredictable points and trajectories}
 
In this section, we will introduce unpredictable points and mention some properties of the corresponding motions. The results will be provided for semi-flows on $X,$ but they are valid for flows as well. We will denote by $\mathbb N$ the set of natural numbers. 
 
\begin{definition}  \label{SSP_point}
A point $p \in X$ and the trajectory through it are \textit{unpredictable} if there exist a positive number $\epsilon_0$ (the sensitivity constant) and sequences $\left\{t_n\right\}$  and  $\left\{\tau_n\right\},$ both of which diverge to infinity, such that
  $\displaystyle \lim_{n \to \infty} f(t_n,p)=p$  and
  $d[f(t_n+\tau_n,p), f(\tau_n,p)] \ge \epsilon_0$ for each $n \in \mathbb N.$
\end{definition} 
 
An important point to discuss is the sensitivity or unpredictability. In the famous research studies \cite{poincare,li,lorenz,smale,poincare1,Devaney}, sensitivity has been considered as a property of a system on a certain set of initial data since it compares the behavior of at least couples of solutions. The above definition allows to formulate unpredictability for a single trajectory.  Indicating   an unpredictable point $p,$ one can make an error by taking a point $f(t_n, p).$ Then $d[f(\tau_n,f(t_n, p)), f(\tau_n,p)] \ge \epsilon_0,$ and this is unpredictability for the motion. Thus,  
we   say   about   the unpredictability   of a single   trajectory  whereas the  former definitions  considered   the  property  in a set  of motions.  In Section \ref{PSC_sec3},  it   will be shown how to  extend  the  unpredictability   to   a chaos.
  
The following assertion is valid. 
  
\begin{lemma} \label{lemma1}
If $p\in X$ is an unpredictable point, then $\mathcal{T}^+(p)$ is neither a rest point nor a cycle.
\end{lemma}   

\noindent \textbf{Proof.} 
Let the number $\epsilon_0$ and the sequences $\left\{t_n\right\},$   $\left\{\tau_n\right\}$ be as in Definition \ref{SSP_point}. Assume that there exists a positive number $\omega$ such that $f(t+\omega,p)= f(t,p)$ for all $t\in \mathbb T_+.$ According to the continuity of $f(t,p),$ there exists a positive number $\delta$ such that if $d[p,q]<\delta$ and $0 \le t  \le \omega,$ then $d[f(t,p), f(t,q)]<\epsilon_0.$ Fix a natural number $n$ such that $d[p_n,p]<\delta,$ where $p_n= f(t_n,p).$ One can find an integer $m$ and a number $\omega_0$ satisfying $0\le \omega_0<\omega$ such that $\tau_n=m\omega+\omega_0.$ In this case, we have that
\[
d[f(\tau_n,p_n),f(\tau_n,p)]=d[f(\omega_0,p_n),f(\omega_0,p)]<\epsilon_0.
\]
But, this is a contradiction since
\[
d[f(\tau_n,p_n),f(\tau_n,p)]=d[f(t_n+\tau_n,p),f(\tau_n,p)] \ge \epsilon_0.
\]
Consequently, $\mathcal{T}^+(p)$ is neither a rest point nor a cycle.
$\square$ 
  
It is seen from the next lemma that the unpredictability can be transmitted  by  the flow. 
  
\begin{lemma} \label{lemma2}
If a point $p\in X$ is unpredictable, then every point of the trajectory $\mathcal{T}^+(p)$ is also unpredictable.
\end{lemma} 
 
\noindent \textbf{Proof.} 
Suppose that the number $\epsilon_0$ and the sequences $\left\{t_n\right\},$   $\left\{\tau_n\right\}$ are as in Definition \ref{SSP_point}. Fix an arbitrary point $q \in \mathcal{T}^+(p)$ such that $q = f(\overline{t},p)$ for some $\overline{t} \in \mathbb T_+.$ One can verify that 
\[
\lim_{n \to \infty} f(t_n,q) = \lim_{n \to \infty} f(t_n + \overline{t},p) =  \lim_{n \to \infty} f(\overline{t}, f(t_n,p)) = f(\overline{t},p)=q.
\]
Now, take a natural number $n_0$ such that $\tau_n>\overline{t}$ for each $n\ge n_0.$ If we denote $\zeta_n = \tau_{n} - \overline{t},$ then we have for $n \ge n_0$ that
\begin{eqnarray*}
& d[f(t_n+\zeta_n,q), f(\zeta_n,q)] & =  d[f(t_n+\zeta_n,f(\overline{t},p)), f(\zeta_n,f(\overline{t},p))] \\
&& = d[f(t_n+\tau_n,p), f(\tau_n,p)] \\
&& \ge \epsilon_0.
\end{eqnarray*}
Clearly, $\zeta_n \to \infty$ as $n \to \infty.$ Consequently, the point $q$ is unpredictable.
$\square$

\begin{remark} 
It is worth noting that the sensitivity constant $\epsilon_0$ is common for each point on an unpredictable trajectory. 
\end{remark}

\section{Chaos on the quasi-minimal set} \label{PSC_sec3}

This section is devoted to the demonstration of chaotic dynamics on a quasi-minimal set. According to \cite{lorenz,Devaney}, the dynamics on a set $S \subseteq X$ is sensitive if there exists a positive number $\epsilon_0$ such that for each $u \in S$ and each positive number $\delta$ there exist a point $u_{\delta} \in S$ and a positive number $\tau_{\delta}$ such that $d[u_{\delta},u]<\delta$ and $d[f(\tau_{\delta},u_{\delta}),f(\tau_{\delta},u)]\ge \epsilon_0.$

The main result of the present study is mentioned in the next theorem, and it is valid for both flows and semi-flows on $X.$

\begin{theorem} \label{thm2}
The dynamics on $\Omega_p^+$ is sensitive if $p \in X$ is an unpredictable point.
\end{theorem}

\noindent \textbf{Proof.}
Let $\epsilon_0 > 0$ be the sensitivity constant corresponding to the point $p.$ Fix an arbitrary positive number $\delta,$ and take a point $r \in \Omega_p^+.$  First of all, consider the case $r \in \mathcal{T}^+(p).$ By Lemma \ref{lemma2}, there exist sequences $\left\{t_n\right\}$ and $\left\{\tau_n\right\},$ both of which diverge to infinity, such that $\lim_{n \to \infty} f(t_n,r)=r$ and $d[f(t_n+\tau_n,r), f(\tau_n,r)] \ge \epsilon_0$ for each $n.$ Fix a natural number $n_0$ such that $d[\overline{r},r]<\delta,$ where $\overline{r}=f(t_{n_0},r).$ In this case, the inequality $d[f(\tau_{n_0},\overline{r}), f(\tau_{n_0},r)] \ge \epsilon_0$ is valid.

On the other hand, suppose that $r \in \Omega^+_p/\mathcal{T}^+(p).$ One can find a sequence $\left\{\eta_m\right\},$ $\eta_m \to \infty$ as $m \to \infty,$ such that $\lim_{m \to \infty} r_m=r,$ where $r_m=f(\eta_m,p).$ According to Lemma \ref{lemma2}, for each $m \in \mathbb N,$ there exist sequences $\left\{s^m_n\right\}$ and $\left\{\xi^m_n\right\},$ both of which diverge to infinity, such that   $\lim_{n \to \infty} r^m_n = r_m$ and $d[f(\xi^m_n,r^m_n),f(\xi^m_n,r_m)] \ge \epsilon_0,$ $n \in \mathbb N,$ where $r^m_n = f(s^m_n,r_m).$

Now, let $m_0$ be a natural number such that $d[r_{m_0},r]<\delta/2.$ Suppose that there exists a natural number $n_1$ satisfying $$d[f(\xi_{n_1}^{m_0},r_{m_0}),f(\xi^{m_0}_{n_1},r)] \ge \epsilon_0/2.$$ If this is the case, then sensitivity is proved. Otherwise, fix $n_2 \in \mathbb N$ such that $d[r^{m_0}_{n_2},r_{m_0}]< \delta/2$ so that $d[r^{m_0}_{n_2},r] \le d[r^{m_0}_{n_2},r_{m_0}]+d[r_{m_0},r] < \delta.$
One can confirm that
\begin{eqnarray*} 
&& d[f(\xi^{m_0}_{n_2},r^{m_0}_{n_2}),f(\xi^{m_0}_{n_2},r)] \\ 
&& \ge  d[f(\xi^{m_0}_{n_2},r^{m_0}_{n_2}),f(\xi^{m_0}_{n_2},r_{m_0})] - d[f(\xi^{m_0}_{n_2},r_{m_0}),f(\xi_{n_2}^{m_0},r)] \\
&& > \epsilon_0/2.
\end{eqnarray*} 

The theorem is proved.
$\square$

In Theorem \ref{thm2}, we have proved the presence of sensitivity in the set $\Omega^+_p$ if $p$ is an unpredictable point in $X.$ In the case that $f$ is a flow on $X,$ one can use the same proof for the verification of sensitivity in $\Omega_p.$ According to Theorem \ref{thm1a} and Lemma \ref{lemma1}, if the positive semi-trajectory of an unpredictable point $p \in X$ is contained in a compact subset of $X,$ then $\Omega^+_p$ contains an uncountable set of everywhere dense stable $P^+$ motions. Additionally, since $\mathcal{T}^+(p)$ is dense in $\Omega^+_p,$ the transitivity is also valid in the dynamics.  

In the definition of Devaney chaos \cite{Devaney}, periodic motions constitute a dense subset. However, in our case, instead of periodic motions, Poisson stable motions take place in the dynamics. Therefore, summarizing the last discussions, we propose a new definition of chaos  based on the concept  of the   unpredictable point. More precisely, we say that the dynamics on the quasi-minimal set $\Omega_p$ is chaotic if the dynamics on it is sensitive, transitive and there exists a continuum of Poisson stable trajectories dense in the quasi-minimal set. Nevertheless, in the framework of chaos there may be infinitely many periodic motions. For instance, chaos in the sense of Devaney \cite{Devaney} and Li-Yorke \cite{li} admit a basis consisting of periodic motions. However, our definition does not contradict to this, and this possibility is exemplified in the next section.

\section{Applications}

In this section, we will mainly investigate symbolic dynamics \cite{Devaney,Wiggins88} and show the presence of unpredictable points as well as chaos on a quasi-minimal set in the sense mentioned in Section \ref{PSC_sec3}. Moreover, we will reveal by means of the topological conjugacy that the same is true for the logistic, H\'enon and horseshoe maps.

Let us take into account the following space of bi-infinite sequences \cite{Wiggins88},
	$$\Sigma^2 = \{s = (\ldots s_{-2} s_{-1}.s_0s_1s_2\ldots): s_j = 0\, \textrm{or}\, 1\ \textrm{for each}\ j\}$$
with the metric 
\[
d[s,\overline{s}] = \sum_{k=-\infty}^{\infty} \frac{|s_k-\bar s_k|}{2^{\left|k\right|}},
\]
where $s = (\ldots s_{-2} s_{-1}.s_0s_1s_2\ldots),$ $\overline{s}=  (\ldots \overline{s}_{-2} \overline{s}_{-1} . \overline{s}_0 \overline{s}_1 \overline{s}_2 \ldots) \in \Sigma^2.$ The shift map $\sigma: \Sigma^2 \to \Sigma^2$ is defined as $$\sigma(\ldots s_{-2} s_{-1}.s_0s_1s_2\ldots) = (\ldots s_{-2}  s_{-1}s_0.s_1s_2\ldots).$$ The map $\sigma$ is continuous and the metric space $\Sigma^2$ is compact \cite{Wiggins88}.

In order to show that the map $\sigma$ possesses an unpredictable point in $\Sigma^2,$ we need an ordering on the collection of finite sequences of $0$'s and $1$'s as follows \cite{Wiggins88}. Suppose that two finite sequences $s = \{s_1 s_2\ldots s_k\}$ and $\bar s = \{\bar s_1 \bar s_2 \ldots \bar s_{k'}\}$ are given. If $k < k',$ then we say that $s < \bar s.$ Moreover, if $k=k',$ then $s < \bar s$ provided that $s_i < \bar s_i,$ where $i$ is the first integer such that $s_i \not = \bar s_i.$ Note that there are $2^m$ distinct sequences of $0$'s and $1$'s with length $m.$ Thus, one can denote the sequences having length $m$ as $s_1^m<\ldots < s_{2^m}^m,$ where the superscript represents the length of the sequence and the subscript refers to a particular sequence of length $m$ which is uniquely specified by the above ordering scheme.
 
Now, consider the following sequence,
\[
s^* = (\ldots s_8^3s_6^3s_4^3s_2^3s_4^2s_2^2.s_1^1s_1^2s_3^2s_1^3s_3^3s_5^3s_7^3\ldots).
\] 
It was demonstrated in \cite{Wiggins88} that the trajectory of $s^*$ is dense in $\Sigma^2.$ We will show that $s^*$ is an unpredictable point of the dynamics $(\Sigma^2, \sigma).$ For each $n\in \mathbb N,$ one can find $j \in \mathbb N$ such that 
$$
s^{2n+2}_{2j-1} = (s^*_{-n} \ldots s^*_0 \ldots s^*_n 0) 
$$ 
and 
$$
s^{2n+2}_{2j} = (s^*_{-n} \ldots s^*_0 \ldots s^*_n 1).
$$ 
Therefore, there exists a sequence $\left\{t_n\right\}$ with $t_n \ge n+ \sum_{k=1}^{2n+1} k 2^{k-1},$ $n \in \mathbb N,$ such that $\displaystyle  \sigma^{t_n}(s^*)=s_i^*$ for $\left|i\right| \le n.$ Accordingly, the inequality $\displaystyle d[\sigma^{t_n}(s^*),s^*] \le 1/2^{n-1}$ is valid so that $\sigma^{t_n}(s^*) \to s^*$ as $n \to \infty.$ Hence, $s^*$ is stable $P^+.$ In a similar way, one can confirm that $s^*$ is stable $P^-.$ Note that $\Sigma^2$ is a quasi-minimal set since $s^*$ is Poisson stable. On the other hand, suppose that there exists a natural number $n$ such that $$\sigma^{t_n+n+1}(s^*)_i = \sigma^{n+1}(s^*)_i$$ for each $i\ge 0.$ Under this assumption we have that $\sigma^{t_n}(s^*)_i=s^*_i$ for $i \ge -n.$ This is a contradiction since the sequence $s^*$ is not eventually periodic. For this reason, for each $n \in \mathbb N,$ there exists an integer $\tau_n \ge n+1$ such that $\sigma^{t_n+\tau_n}(s^*)_0 \neq \sigma^{\tau_n}(s^*)_0.$ Hence, $d[\sigma^{t_n+\tau_n}(s^*),\sigma^{\tau_n}(s^*)]\ge 1$ for each $n \in \mathbb N,$ and $s^*$ is an unpredictable point in $\Sigma^2.$  

One of the concepts that has a great importance in the theory of dynamical systems is the topological conjugacy, which allows us to make interpretation about complicated dynamics by using simpler ones. Let $X$ and $Y$ be metric spaces. A flow (semi-flow) $f$ on $X$ is topologically conjugate to a flow (semi-flow) $g$ on $Y$ if there exists a homeomorphism $h:X \to Y$ such that $h \circ f = g \circ h$ \cite{Devaney,Wiggins88}. The following theorem can be verified by using the arguments presented in \cite{Banks}.

\begin{theorem} \label{PSC_thm3}
Suppose that $X$ and $Y$ are metric spaces and a flow (semi-flow) $f$ on $X$ is topologically conjugate to a flow (semi-flow) $g$ on $Y.$ If there exists an unpredictable point whose trajectory is contained in a compact subset of $X,$ then there also exists an unpredictable point whose trajectory is contained in a compact subset of $Y.$
\end{theorem}

Since the shift map $\sigma$ on $\Sigma^2$ is topologically conjugate to the Smale Horseshoe \cite{Devaney,Wiggins88}, one can conclude by using Theorem \ref{PSC_thm3} that the horseshoe map possesses an unpredictable point and a trajectory. On the other hand, let us consider the H\'enon map
\begin{eqnarray} \label{Henon_map}
\begin{array}{l}
x_{n+1}=\alpha - \beta y_n-x_n^2 \\
y_{n+1} = x_n,
\end{array}
\end{eqnarray}
where $\beta \neq 0$ and $\alpha \ge (5+2\sqrt{5}) (1+\left|\beta\right|)^2 /4.$ 
It was proved by Devaney and Nitecki \cite{Dev79} that the map (\ref{Henon_map}) possesses a Cantor set in which the map is topologically conjugate to the shift map $\sigma$ on $\Sigma^2.$ Therefore, Theorem \ref{PSC_thm3} also implies the  presence of  an unpredictable point and a trajectory in the dynamics of (\ref{Henon_map}).

%%%%%%%%

Next, as an example of a semi-flow, consider the following space of infinite sequences \cite{Devaney},
\[
\Sigma_2 = \{s = (s_0s_1s_2\ldots): s_j = 0\, \textrm{or}\, 1\ \textrm{for each}\ j\}
\]
with the metric 
\[
d[s,\overline{s}] = \sum\limits_{k=0}^{\infty}\frac{|s_k-\overline{s}_k|}{2^k},
\]
where $s=(s_0 s_1 s_2\ldots),$ $\overline{s}=(\overline{s}_0 \overline{s}_1 \overline{s}_2\ldots) \in \Sigma_2.$ The shift map $\sigma: \Sigma_2 \to \Sigma_2$ is defined as
$\sigma(s_0s_1s_2\ldots) = (s_1s_2s_3\ldots).$ As in the case of the space of bi-infinite sequences, the metric space $\Sigma_2$ is compact and the map $\sigma$ is continuous \cite{Devaney,Wiggins88}.

Let us take into account the sequence 
\[
s^* = (\underbrace{0 \ 1}_{1 \, blocks}|\underbrace{00 \ 01 \ 10 \ 11}_{2 \,  blocks}|\underbrace{000 \ 001 \ 010 \ 011 \ \ldots}_{3 \, blocks}| \ldots),
\]
which is constructed by successively listing all blocks of $0$'s and $1$'s of length $n,$ then length $n+1,$ etc. This sequence is non-periodic and its trajectory $\mathcal{T}(s^*)=\left\{\sigma^{n}(s^*): \ n=0,1,2,\ldots \right\}$ is dense in $\Sigma_2$ \cite{Devaney}. Note that the number of all blocks of length $n$ in $s^*$ is $2^n.$ Based upon the construction of $s^*,$ there exists a sequence $\left\{t_n\right\}$ satisfying $t_n \ge \sum_{j=1}^n j2^j,$ $n \in \mathbb N,$ such that $s^*_i=\sigma^{t_n}(s^*)_i$ for each $i=0,1,2,\ldots,n.$ Clearly, $t_n \to \infty$ as $n \to \infty$ and $d[\sigma^{t_n}(s^*),s^*] \le 1/2^{n}$ so that $\sigma^{t_n}(s^*) \to s^*$ as $n \to \infty.$ Hence, $s^*$ is stable $P^+.$ In a very similar way to the bi-infinite sequences, one can show the existence of a sequence $\left\{\tau_n\right\},$ $\tau_n \to \infty$ as $n \to \infty,$ such that $d[\sigma^{t_n+\tau_n}(s^*),\sigma^{\tau_n}(s^*)] \ge 1$ for each $n \in \mathbb N.$ Thus, $s^*$ is an unpredictable point in $\Sigma_2.$ 

It was shown in \cite{rob} that the logistic map $x_{n+1} = \mu x_n (1-x_n)$ possesses an invariant Cantor set $\Lambda \subset [0,1],$ and the map on $\Lambda$ is topologically conjugate to $\sigma$ on $\Sigma_2$ for $\mu >4.$ Therefore, the map with $\mu >4$ possesses an unpredictable point and a trajectory in accordance with Theorem \ref{PSC_thm3}.

\section{Conclusions}  

Emphasizing the ingredients of Devaney \cite{Devaney} and Li-Yorke \cite{li} chaos, one can see that the definitions of chaos have been considered by means of sets of motions, but not through a single motion description. This is the first time in the literature that we initiate chaos by a single function. Thus, the line, equilibrium, periodic function, quasi-periodic function, almost periodic function, recurrent function, Poisson stable motion is prolonged with the new element - unpredictable motion. This supplement to the line creates the possibility of other functions behind the known ones.

An essential point to discuss is the sensitivity or unpredictability. In the literature, sensitivity has been considered through initially nearby different motions. However, we say about unpredictability as an interior property of a single trajectory. Then chaos appears in a neighborhood of the trajectory. The symbolic dynamics illustrates all the results and it is an important tool for the investigation of the complicated dynamics of continuous-time systems such as the Lorenz and R\"{o}ssler equations \cite{Mischaikow,Zgliczynski}. It is worth noting that unpredictable points can be replicated by the techniques summarized in \cite{Akh1}.

One can see the proximity of chaos and quasi-minimal sets by comparing their definitions \cite{li,Nemytskii,Devaney}. Transitivity is a common feature of them, and the closure of a Poisson stable  trajectory contains infinitely many Poisson stable orbits. In its own turn we know that a periodic trajectory is also Poisson stable. Possibly, the uncountable set of Poisson stable trajectories is an  ultimate form of infinitely many cycles known for a chaos. One may also ask whether sensitivity is proper for any quasi-minimal set. These are the questions to discover more relations between quasi-minimal sets and chaos.

\end{document}